\documentclass{article}%
\usepackage{latexsym,amssymb,amsmath,amscd,amsthm,amsxtra}
\usepackage{amsmath}
\usepackage{amsfonts}
\usepackage{amssymb}
\usepackage{pb-diagram }
\usepackage{pb-xy}
\usepackage[all]{xy}
\usepackage{graphicx}%
\newtheorem{thm}{Theorem}[section]
\newtheorem{lem}[thm]{Lemma}
\newtheorem{prop}[thm]{Proposition}
\newtheorem{co}[thm]{Corollary}
\newtheorem{defi}[thm]{Definition}

\newtheorem{rk}[thm]{Remark}
\newtheorem{ex}[thm]{Example}

\newcommand{\zz}{\mathbb {Z}}
\newcommand{\reals}{\mathbb {R}}
\newcommand{\integers}{\mathbb {Z}}
\newcommand{\torus}{{\mathbb T}}
\newcommand{\cala}{{\mathcal A}}
\newcommand{\cald}{{\mathcal D}}
\newcommand{\calt}{{\mathcal T}}
\newcommand {\pf}{\noindent{\bf Proof.}\ }
\newcommand{\double}{\reals^n\oplus {\reals^n}^*}
\newcommand{\arrows}{\,\lower1pt\hbox{$\longrightarrow$}\hskip-.24in\raise2pt
             \hbox{$\longrightarrow$}\,}

\begin{document}

\title{Quantization and Morita Equivalence for
 Constant Dirac Structures on Tori}
\author{Xiang Tang and Alan Weinstein\thanks{Research of both
authors partially supported
by NSF Grant DMS-02-04100.
\newline \mbox{}\ \ \ \  MSC2000 Subject Classification Numbers: 46L65,81S10
(Primary). 
\newline \mbox{}\ \ \ \   Keywords: Dirac structure, Poisson structure,
Morita equivalence, quantization.}
\\Department of Mathematics\\ University of California\\ Berkeley, CA
94720 USA\\ {\small(xtang@math.berkeley.edu, alanw@math.berkeley.edu)}}
\date{}
\maketitle

\begin{abstract}
We define a $C^*$-algebraic quantization of constant Dirac structures on
tori, which
extends the standard quantization of Poisson structures.  We prove
that Dirac structures in the same orbit of a natural action of
$O(n,n|\mathbb{Z})$ give rise to Morita equivalent algebras,
completing and extending a theorem of Rieffel and Schwarz.
\end{abstract}

\section{Introduction}

Quantum, or noncommutative, tori can be obtained by deformation
quantization of constant Poisson structures on ordinary tori.  This
fact was noticed by the second author \cite{we:rotation} and was
developed by Rieffel \cite{ri:deformation} into a rigorous theory of
``strict deformation quantization'' from Poisson manifolds to $C^{*}$
algebras.  Earlier, Rieffel \cite{ri:morita1} had introduced a notion
of strong Morita equivalence for operator algebras, sufficient to
imply the equivalence of suitable categories of topological
representations.  For simplicity, we shall refer to this notion simply
as ``Morita equivalence.''

Motivated by applications to duality in string theory (see
\cite{sc:morita}), Rieffel and Schwarz \cite{ri-sc:morita} showed,
with an additional technical hypothesis, that the algebras of
functions on two noncommutative $n$-tori are Morita equivalent if the
underlying Poisson structures are related by a ``fractional linear
transformation'' whose coefficient matrix belongs to
$SO(n,n|\mathbb{Z}).$ In this paper, we will prove the Rieffel-Schwarz
result without the additional hypothesis by extending the scope of the
theorem from Poisson structures to Dirac structures, whose definition
we will recall later in this introduction.

Li \cite{li:strong} has also proven the full theorem for the case of
Poisson structures, and in fact all of our results could be deduced from
his theorem.  Nevertheless, our proof is completely different, and we
believe that it sheds further light on the role of
$SO(n,n|\mathbb{Z}).$

A constant Poisson structure on $\torus^n=\reals^n/\integers^n$ is
specified by a skew-symmetric, real, $n\times n$ matrix $\Pi$,
representing a skew-symmetric bilinear form on the space
${\reals^n}^*$ of translation-invariant vector fields on the torus.
The standard deformation quantization of the Poisson manifold
$(\torus^n,\Pi)$ is obtained by deforming the algebra of finite
Fourier series on $\torus^n$ using the multiplication rule
\begin{equation}
\label{eq-moyal}
e_r *_{\hbar} e_s = e^{-\pi i\hbar \Pi(r,s)}e_{r+s},
\end{equation}
where $r$ and $s$ are multi-indices in $\integers^n,$ $\hbar$ is the
deformation parameter, and
$e_m(x) = e^{2\pi imx}$.
This product extends to the Fourier series with rapidly decreasing
coefficients and from there to a $C^*$ completion
$\cala_{\hbar\Pi}$ which is known as ``the algebra of continuous functions
on
the
quantum torus\footnote{Sometimes $\cala_{\hbar\Pi}$ is  itself referred
  to as the quantum torus, but since this terminology does not agree
  with normal usage when $\Pi=0$, we prefer not to use it.}
$\torus_{\hbar\Pi}$''.  At this point, we will leave the world of
deformations by setting $\hbar=1$ so that we have, for each
skew-symmetric form $\Pi$, the algebra $\cala_\Pi$ of functions on
$\torus_\Pi$.  This algebra may also be described as the $C^*$-algebra
determined by $n$ unitary generators $\epsilon_1,\ldots,\epsilon_n$ (lower
indices here are in  $\integers$ rather than $\integers^n$) subject to
the commutation relations $\epsilon_j \epsilon_i = e^{-2\pi i\hbar
\Pi_{ij}} \epsilon_i \epsilon_j.$

One sees immediately from the commutation relations that
two kinds of operations on $\Pi$.
do not change the isomorphism class of $\cala_\Pi$.
Adding a matrix with integer
coefficients does not change the algebra at all.  Also, if
we $A \in GL(n,\integers)$, then the algebras $\cala_\Pi$ and
$\cala_{A\Pi A^t}$ are
isomorphic via the map which takes $\epsilon_i$ to
$\prod_j \epsilon_j^{a_{ij}}.$

Much less obvious are operations for which the quantized algebras are
not isomorphic but are still Morita equivalent.  The first of these was
discovered by Rieffel \cite{ri:irrational} when  $n=2$.
The matrix
$\Pi$ then has the form
$\left(
\begin{array}
[c]{cc}
0 & \theta\\
-\theta & 0
\end{array}  \right) , $
where $\theta$ is a real number.  The first of the two kinds of
operation
above adds an
integer to $\theta$, while the second leaves $\theta$ fixed or simply
changes its sign.  Rieffel proved in \cite{ri:irrational} that
a third operation, namely replacing
$\theta$ by $1/\theta$,  preserves the Morita equivalence
class of the algebra.   The three types of
operations are contained in, and in fact generate, the action of the
group $GL(2,\integers)$ on the real numbers by
fractional linear transformations $\theta\mapsto (p\theta + q)(r\theta
+s)^{-1};$ hence, the algebras
corresponding to $\theta$ and $\theta'$  are Morita equivalent if
$\theta$ and $\theta'$ are in the same $GL(2,\integers)$ orbit.  The
converse is also proved in  \cite{ri:irrational}.

To a large extent, the main result of \cite{ri-sc:morita} extends
to higher dimensional tori the
 ``if'' part of the classification above
(but not the converse, which is false unless $n=2$).
The group $O(n,n|\reals)$ of automorphisms of $\reals^n \oplus {\reals
 ^n}^*$ preserving  the
indefinite inner product
\begin{equation}
\label{eq-bilinear}
(X_{1} + \xi_{1}, X_{2} +\xi_{2})=\frac{1}{2}(<\xi_{1}, X
_{2}>+ <\xi_{2}, X _{1}>),
\end{equation}
has subgroups $O(n,n|\integers)$ and
$SO(n,n|\integers)$ defined in the obvious way.
An element of $O(n,n|\integers)$ may be written in block form as
$$
g = \left( \begin{array}[c]{ll}
A & B\\
C & D
\end{array} \right),
$$
where $A$, $B$, $C$, and $D$ are $n\times n$ integer matrices which
satisfy $A^tC+C^tA=0=B^tD+D^tB$ and $A^tD+C^tB=1$.
Such a matrix  ``acts''  on the space ${\calt}_n$ of skew-symmetric
$n\times n$
matrices
by taking $\Pi$ to
\begin{equation}
\label{eq-action}
g\cdot\Pi= (A\Pi + B)(C\Pi +D)^{-1}  .
\end{equation}
The word ``acts'' is in quotation marks because the right hand side is not
defined when the denominator $C\Pi + D$ is singular.

The
main result of \cite{ri-sc:morita} may now be stated as follows.

\begin{thm}
\label{thm-risc}
If $\Pi\in {\calt}_n$ is such that $g\cdot\Pi$ is defined for all
$g\in SO(n,n|\integers),$ then the algebras $\cala_{g\cdot\Pi}$ are all
Morita equivalent to one another.
\end{thm}

We note that \cite{ri-sc:morita} contains a counterexample to the
converse of this theorem when $n=3$, while Schwarz \cite{sc:morita}
proves a converse using a refined notion known as ``complete Morita
equivalence.'' 

The proof  of Theorem \ref{thm-risc} in \cite{ri-sc:morita}
 uses a decomposition of
the general element $g\in SO(n,n|\integers)$ as a product of generators of
three types, analogous to the three types described above for $n=2$.
  As a result, it does not establish the Morita
equivalence of $\Pi$ and $g_0\cdot\Pi$ if the action on $\Pi$ is defined
for a
particular $g_0$ but not for all $g$.  The key idea of the present
paper is to circumvent this difficulty by enlarging (in fact,
compactifying)
${\calt}_n$ to the space of Dirac structures, on
which the action of $ O(n,n|\reals)$, and hence that of
$SO(n,n|\integers),$ is everywhere defined.   This idea was suggested by
the appearance of the bilinear form (\ref{eq-bilinear}) in both the
Rieffel-Schwarz theorem and the definition by
Courant \cite{co:dirac} of Dirac structures, which we now recall.

\begin{defi}
\label{def-dirac}
A {\bf Dirac structure} on a vector space $V$ is a maximal
isotropic subspace of $V\oplus V^*$ with respect to the
non-degenerate symmetric bilinear form
\eqref{eq-bilinear}.
  A Dirac structure on a vector bundle $E$ is a
subbundle of $E\oplus E^*$ which is a Dirac structure on each fibre.
A Dirac structure on a manifold $M$ is a Dirac structure on $TM$ whose
space of sections is closed under the {\bf Courant bracket}
\begin{equation}
\label{eq-courant}
[X_1+ \xi_1, X_2 +\xi_2]=
[X_1,X_2]+ L_{X_1}\xi_2-L_{X_2}\xi_1 +
{\textstyle{\frac 12}}d(<\xi_1,X_2>-<\xi_2,X_1>).
\end{equation}
The space of all Dirac structures on the vector space $\reals^n$ will be
denoted by $\cald_n$.
\end{defi}

Dirac structures on a manifold $M$ include the Poisson
structures and closed 2-forms (identified with the graphs of the
corresponding skew-symmetric bundle maps $T^*M \to TM$ or $TM \to
T^*M$), as well as the foliations (identified with direct sums $F\oplus
F^\circ,$ where $F$ is an integrable subbundle of $TM$, and $F^\circ$
is its annihilator).

  From now on, we will be concerned exclusively with
constant (i.e. translation-invariant) Dirac structures on tori
$\torus^n$.  Since all the terms in the Courant bracket involve
derivatives, these are just the translation-invariant Dirac structures
on the tangent bundle $T\torus^n$, or equivalently the Dirac
structures on the vector space $\reals^n$ of constant vector fields
on $\torus^n$ (which
may be identified with the tangent space at any point).

The ``action'' of $O(n,n|\reals)$ on antisymmetric matrices now has the
following interpretation.  We identify each $\Pi \in \calt_n$ with
the map ${\reals^n}^* \to \reals^n$ of which it is the matrix with respect
to
the standard basis and its dual.  Then we
identify   $\Pi$ with its graph
$$\Gamma_\Pi = \{(\Pi \xi,\xi)|\xi\in {\reals^n}^*\},$$
an element of $\cald_n$.  The group $O(n,n|\reals)$ acts in the obvious
way
on $\cald_n$, and the correspondence $\Pi\mapsto \Gamma_\Pi$ is
$O(n,n|\reals)$ equivariant with respect to the action
\eqref{eq-action}.
Thus, if a product of generators $g=g_r\cdots g_1$ maps
$\Pi$ to $\Pi'$, even if $g_1\Pi$,
$g_2g_1\Pi$,\ldots are not all defined as antisymmetric matrices, they
{\em are} defined as Dirac structures.  Our strategy for proving the
Rieffel-Schwarz theorem, then, is to attach an algebra $\cala_\Gamma$
 (more precisely, a
Morita equivalence classes of algebras) to each
Dirac structure
$\Gamma$, and to prove that this Morita equivalence class is
unchanged when $\Gamma$ is transformed by any
member of a certain set of generators of
$O(n,n|\integers)$.

\begin{rk}
\label{rk-oso}
{\em
The apparent extension of the Rieffel-Schwarz theorem from
$SO(n,n|\integers)$ to $O(n,n|\integers)$  is  illusory.  As we
will see in Corollary \ref{co-parity} below
any $g \in O(n,n|\integers)$ which transforms some Poisson
structure into another one must lie in $SO(n,n|\integers)$.
  On the other hand, by passing to $O(n,n|\integers)$, we will not only
  bypass the ``obstruction'' in the original proof, but we will reduce
  from three to two the number of kinds of generators which must be dealt
with.  }
\end{rk}

To  construct $\cala_\Gamma$, we begin by recalling
   from \cite{co:dirac} (Proposition 1.1.4) that to every
Dirac structure $\Gamma$ on the vector space $\reals^n$ there corresponds
a natural bivector $\Pi_{\Gamma}$ on the quotient
of $\reals^{n}$ by $\reals^{n}\cap\Gamma$.  (Here, we are
identifying
$\reals^n$ with the subspace $\reals^n\oplus \{0\}$ of
$\reals^n\oplus {\reals^n}^*.$)  Conversely (with a proof
similar to arguments in Section 1.1 of \cite{co:dirac}), the intersection
$\reals^{n}\cap\Gamma$ and the bivector $\Pi_{\Gamma}$ determine
$\Gamma$.

For a constant Dirac structure $\Gamma$ on $\torus^n$, the intersection
$\reals^{n}\cap\Gamma$ defines a foliation known as the
characteristic foliation of the Dirac structure, while the bivector
$\Pi_\Gamma$ defines a Poisson structure transverse to this
foliation.  The  idea behind what follows is that the function algebra
of the quantized Dirac manifold $(\torus^n,\Gamma)$
should be obtained from a foliation algebra for
the characteristic foliation (specifically, the groupoid algebra of the
leafwise fundamental groupoid of the foliation) by ``deformation''
via the transverse Poisson
structure.  This idea, already suggested by Block and
Getzler \cite{bl-ge:quantization} (see also \cite{xu:noncommutative}),
is at best inconvenient to apply directly in our setting, since it
requires the extra data of a so-called Haar form.  Instead, we restrict
the fundamental groupoid of the foliation
to a subtorus $M\subseteq \torus^n$ which
is a complete transversal to obtain an equivalent
\'etale
groupoid.  $M$ carries a Poisson structure
$\Pi_{\Gamma,M}$, from
which we can obtain a quantum torus
 algebra $\cala_{\Pi_{\Gamma,M}}$ in the usual way.
The restricted groupoid turns out to be a transformation groupoid for
an action of a lattice  on $M$, and so we may form the
crossed product algebra of $\cala_{\Pi_{\Gamma,M}}$ with the lattice.
This
algebra depends on the choice of $M$, but we will show
 that its Morita equivalence class depends only on $\Gamma$.  This
 ``independence of transversal'' proof, carried out in detail (in a
 more general context) in \cite{ta:lemma}, is the only analytic
 ingredient in our proof.  Its use is thus an extension of the idea
of Connes (Section $8.\beta$ in \cite{C}), that the original
 Rieffel theorem for $\torus^2$ is a consequence of the Morita invariance
 of ordinary foliation groupoid algebras under change of transversal.

{ \textbf{Acknowledgments}: We would like to thank H.~Bursztyn,
M.~Crainic, K.~Fukaya, H.~Li, M.~Rieffel and A.~Schwarz for helpful
suggestions.}

\section{Linear algebra of Dirac structures}
We begin with some notation and definitions related to Dirac
structures on a vector space $V$.

We will identify $V$ with the subspace $V\oplus{0}$ of $V\oplus V^*$.
The projection  $V\oplus V^*\to V^*$ will be denoted by $p_*$.  The
annihilator in $V^*$ of a subspace $E\subseteq V$ will be denoted by
$V^\circ$.

\begin{defi}
The {\bf characteristic subspace} $C(\Gamma)$ of a Dirac structure
$\Gamma$ on $V$ is $\Gamma \cap V$.  The dimension of $C(\Gamma)$ is
the {\bf nullity} $N(\Gamma)$ of $\Gamma$.   $\Gamma$ is called {\bf
  even} or {\bf odd} according to the dimension of its nullity.
\end{defi}

The maximal isotropic property of $\Gamma$ immediately implies that
\begin{equation}
\label{eq-nullproj}
C(\Gamma)^\circ = p_*(\Gamma).
\end{equation}

\begin{lem}
\label{lem-parity}
The spaces of even and odd Dirac structures
are the connected components of
$\cald_n$.  The action of $SO(n,n|\reals)$ leaves each component
invariant, while the action of $O(n,n|\reals)\setminus
SO(n,n|\reals)$ interchanges the components.
\end{lem}

\pf
We follow Courant \cite{co:dirac} and rewrite $\double$ as an orthogonal
direct sum $P\oplus N$, by choosing a positive definite symmetric
inner product $b$ on $\reals^n$ and letting $P$ be the graph of the
corresponding map $\reals^n\to{\reals^n}^*$.  Its orthogonal
complement $N$ is the graph of the map corresponding to $-b$.
$P$ and $N$ are maximal positive definite and negative subspaces
of $\double$, and the Dirac structures on $\reals^n$ are the graphs of
anti-isometries from $P$ to $N$, with which they may be identified.
If we fix one such anti-isometry $K$, namely the one identified with
$\reals^n\oplus \{0\}$, each Dirac structure $\Gamma$ may be
identified by composition with $K^{-1}$ with an
isometry $L$ from $P$ to itself, and the characteristic subspace of
$\Gamma$ is then identified with the fixed space of $L$.  Since
the codimension of the fixed space of $L$ is even or odd
according to whether $L$ preserves or reverses orientation,
the dimension modulo 2 of the characteristic subspace is
constant on each component of $\cald_n$.

To analyze the effect of
$O(n,n|\reals)$ on parity, it suffices to look at the maximal compact
subgroup
$O(n)\times O(n)$ which leaves the subspaces $P$ and $N$ invariant.
An element of this subgroup is a pair $(A,B)$ of isometries of $P$ and
$N$ respectively, and it takes the Dirac structure associated with
$K^{-1}L:P\to P$ to the structure associated with $K^{-1}BLA^{-1}$.
The last part of the lemma
 now follows from the identity $\det K^{-1}BLA^{-1} = \det
K^{-1}L \det B \det A$.
\qed

Since Poisson structures are even Dirac structures, we have the
following corollary, which gives one explanation for the appearance of
the {\em special} orthogonal group in the Rieffel-Schwarz theorem.

\begin{co}
\label{co-parity}
Any element of
$O(n,n|\reals)$ which transforms one Poisson structure to another must
lie in $SO(n,n|\reals)$.
\end{co}

We now introduce further notation in $\reals^n$ and $O(n,n|\integers)$.
Let $(e_1,\ldots,e_n)$ and $(f_1,\ldots f_n)$
be the standard basis of $\reals^n$ and its dual basis.  For any
subset $I$ of $\{1,2,\ldots,n\}$, we denote its complement by $I'$ and
its cardinality by $|I|$.
$\reals^I$ and ${\reals^I}^*$ will denote the subspaces of $\reals^n$
and ${\reals^n}^*$ spanned by the $e_i$ (or $f_i$) for $i\in I$.  Note
that $({\reals^I})^\circ ={\reals^{I'}}^*.$   We will also use this
``exponent'' notation for subgroups  of $\integers^n$ and $\torus^n$.

The
  element of $O(n,n|\integers)$ which exchanges $e_i$ with $f_i$ for
  all $i\in I$ and which leaves all other basis elements fixed will be
  denoted by $\sigma_I$.
If $A$ belongs to $GL(n,\integers)$, we will denote by $\rho(A)$ the
element of $O(n,n|\integers)$ which acts on the first summand
of $\reals^n\oplus {\reals^n}^*$ by $A$ and on the second summand by
$(A^t)^{-1}$.  If $N$ is a skew symmetric $n\times n$ matrix with
integer entries, $\nu(N)$ will denote the map $(x,y)\mapsto
(x+Ny,y)$.  (When applied as a fractional linear transformation,
$\nu(N)$ just adds $N$ to each Poisson structure.)  The additive group
of all $\nu(N)$'s is generated by its elements $\nu_{ij}$ for
$i < j$, where $\nu_{ij}$ is the sum of the identity with the rank 2
matrix which maps $f_i$ to $e_j$, $f_j$ to $-e_i$, and all other basis
elements to zero.

It is proven in \cite{ri-sc:morita} that $SO(n,n|\integers)$ is
generated by the $\rho(A)$'s, the $\nu_{ij}$'s, and $\sigma_{\{1,2\}}$.
It
follows easily that $O(n,n|\integers)$ is generated by
 the $\rho(A)$'s, the $\nu_{ij}$'s, and $\sigma_{\{1\}}$.  But in fact
 even more is true.

\begin{lem}
\label{lem-generators}
The group $O(n,n|\integers)$ is generated by $\sigma_{\{1\}}$ and
$\rho(GL(n,\integers)).$
\end{lem}
\pf
The subgroup generated by $\sigma_{\{1\}}$ and
$\rho(GL(n,\integers))$ contains $\sigma_{\{i\}}$ for all $i$.   Hence
it contains $\sigma_{\{1\}}\sigma_{\{2\}}=\sigma_{\{1,2\}}.$
A straightforward computation shows that $\nu_{ij}=\sigma_{\{i\}} \rho(A)
\sigma_{\{i\}}$, where $A$ maps $e_i$ to $e_i+e_j$ and fixes all the
other basis elements of $\reals^n$.
\qed

\begin{rk}
{\em
Elimination of the generators $\nu_{ij}$ will be essential in our
proof.  Although these generators act trivially on the quantization of
Poisson structures, this is not true for general Dirac structures.
}
\end{rk}

Our quantization of Dirac structures will make essential use of the
following result of elementary linear algebra.

\begin{lem}
\label{lem-subcomplement}
For any $\Gamma \in \cald_n$, there is a subset $I$ with
$|I|=n-C(\Gamma)$ for which $\reals^I$ is complementary to $N(\Gamma)$.
\end{lem}

We will also use:

\begin{lem}
\label{lem-comp}
If $\reals^{I'}$ is complementary to $C(\Gamma)$ (so that, in particular,
 $|I|=N(\Gamma)$),  then
$N(\sigma_I(\Gamma))$=0; i.e. $\sigma_I(\Gamma)$ is a
Poisson structure.
Conversely, if $|I|=N(\Gamma)$ and $\sigma_I(\Gamma)$ is a
Poisson structure, then $\reals^{I'}$ is complementary to
 $C(\Gamma)$.
\end{lem}
\pf
Taking annihilators in the direct sum decomposition
\begin{equation}
\label{eq-directsum}
C(\Gamma)\oplus \reals^{I'} = \reals^n
\end{equation}
and using \eqref{eq-nullproj}  gives the dual decomposition
\begin{equation}
\label{eq-directsumdual}
p_*(\Gamma)\oplus {\reals^I}^* = {\reals^n}^*.
\end{equation}
Now $\sigma_I(\Gamma)$ is a
Poisson structure if it has zero intersection with $\reals^n$ or,
equivalently, if $\Gamma \cap \sigma_I(\reals^n) = \{0\}.$  Suppose,
then, that $(X,\xi) \in \Gamma \cap \sigma_I(\reals^n)$.  Since
$\sigma_I(\reals^n)={\reals}^{I'}\oplus {\reals^I}^*$, it follows from
\eqref{eq-directsumdual} that $\xi=0$.  Now we have $(X,0)\in\Gamma$
with $X\in \reals ^{I'}$, which by \eqref{eq-directsum} implies that
$X=0$ as well.

For the converse, it suffices to prove that, if $\sigma_I(\Gamma)$ is
a Poisson structure, then $\reals^{I'}\cap
C(\Gamma)=\{0\}.$  But if $X$ is an element of this intersection, then
$(X,0)\in \Gamma$, and $\sigma^I(X,0)=(X,0).$  Since
$\sigma_I(\Gamma)$ is a Poisson structure, we must have
$X=0$.
\qed

We will see below that, when the constant Dirac structure  on
$\torus^n$ given by $\Gamma$ is quantized by our method, the resulting
algebra is just
that obtained by quantizing the Poisson structure $\sigma_I(\Gamma)$.

In general, there are many $\sigma_I$'s which can convert a given
Dirac structure to a Poisson structure.  However, the
ones with minimal length are those given by Lemma
\ref{lem-comp}.  This follows from the following lemma and its corollary.

\begin{lem}
\label{lem-nullity}
For any Dirac structure $\Gamma$ and any $i$,
$N(\sigma_{\{i\}}(\Gamma))=N(\Gamma)\pm 1$.
\end{lem}
\pf
It follows from  \eqref{eq-nullproj} that the statement is
true if $\dim p_*(\sigma_{\{i\}}(\Gamma))=\dim
p_*(\Gamma)\pm 1.$  Writing $\sigma_{\{i\}}$ as $1+R$, where $R$ has
rank 2, we find that
$$p_*(\sigma_{\{i\}}(\Gamma))=p_*((1+R)(\Gamma))\subseteq p_*(\Gamma) +
p_* R(\Gamma).$$  Since $p_* R(\Gamma)$ has dimension at most 2, it
follows that $\dim p_*(\sigma_{\{i\}}(\Gamma))\leq \dim
p_*(\Gamma) +2.$  Since $\sigma_{\{i\}}$ is an involution, we also
have $\dim p_*(\sigma_{\{i\}}(\Gamma))\geq \dim
p_*(\Gamma) -2,$ and the result now follows from Lemma \ref{lem-parity}.
\qed

\begin{co}
\label{co-length}
If $\sigma_I(\Gamma)$ is the a Poisson structure, then
  $|I|\geq N(\Gamma)$.
\end{co}

\section{ Quantization}
\label{sec-quantization}
In this section, we will define the quantization of a constant
Dirac structure on a torus as a Morita equivalence class of algebras.
In fact, we will see that quantizations can all be realized as
quantizations of constant {\em Poisson} structures.\footnote{We are
indebted to Marc Rieffel for pointing this out.}

   From now on, we will identify each constant Dirac structure on
$\torus^n$ with the corresponding $\Gamma\in \cald_n$.  The subspace
$C(\Gamma)$ then determines the direction of the characteristic
foliation.

It follows from Lemma \ref{lem-subcomplement} that $C(\Gamma)$ has
a complementary subspace $\widetilde{M}$ having a basis
with rational components;
we will often work with complements of
the special form given by selecting a subset of the coordinate axes,
but any rational complement can be put in this form by a change of
coordinates in $GL(n,\integers)$.  The
projection of any rational complement into $\torus^n$ is a compact
subtorus $M$ which is a complete transversal to the
characteristic foliation.

The Dirac structure $\Gamma$ induces a transverse Poisson structure
$\Pi_{\Gamma,M}$ on the transversal torus $M$.  We will construct an
algebra $\cala_{\Gamma,M}$ by using this Poisson structure to
``quantize'' the groupoid algebra of the
restriction to $M$ of the fundamental groupoid along the leaves of the
characteristic foliation of $\Gamma$.  The full fundamental groupoid
is naturally isomorphic to the transformation groupoid
associated with the translation action  of
$C(\Gamma)$  on $\torus^n$; its restriction to $M$ becomes the
transformation groupoid of the translation action of a lattice.
  The subspace and lattice may be identified with
$\reals^n/\widetilde{M}\cong \reals^k$ and its integer lattice
$\integers^k$ respectively.
 In more geometric terms, we may consider $\torus^n$
 as a principal $M$ bundle over the quotient torus $\torus^n/M$
 whose fundamental group is the lattice $\integers^k$.  The
 characteristic foliation is a flat connection on this principal
 bundle, and the homomorphism $\integers^k \to M$ giving the
 translation action is the  holonomy of this foliation.

The translation action on the torus $M$ clearly induces an action on the
quantum torus algebra $\cala_{\Pi_{\Gamma,M}},$ which enables us to
make the following definition.
\begin{defi}
\label{defi-algebra}
With notation and terminology as above, we define the algebra
$\cala_{\Gamma,M}$ to be
the crossed product $C^*$-algebra
$\cala_{\Pi_{\Gamma,M}}\times\integers^{k}$.
\end{defi}

\begin{rk}
{\em   We
may think of the crossed product  $\cala_{\Gamma,M}
=\cala_{\Pi_{\Gamma,M}}\times\integers^{k}$ as a ``quantization'' of
the algebra $C(M)\times\integers^{k}$.  The latter is the groupoid
algebra of the transformation groupoid $M\times \integers^k\arrows M$
associated to the translation action of $\integers^k$ on $M$.  We may
therefore think of  $\cala_{\Gamma,M}$ as the groupoid algebra (but
not the function algebra) of a quantum groupoid.  It may thus be
considered as a (strict) deformation quantization of a noncommutative
algebra, in the spirit of  \cite{bl-ge:quantization} and
\cite{xu:noncommutative}.  We refer to \cite{ta:lemma} for more details.
}
\end{rk}

The construction just described becomes much more concrete when
the complement $\widetilde{M}$ is taken to be of the form $R^{I'}$
for suitably chosen $I'$.  (As we remarked earlier, any complement can
be put in this form by a transformation in $GL(n,\integers)$.)
$C(\Gamma)$ is then the graph of
a linear map $\beta:\reals^{I}\to \reals^{I'}$ whose composition with
the inclusion $\integers^I\to\reals^I$ and the projection
$\reals^{I'}\to \torus^{I'}$ is the holonomy action of $\integers^I$
on $\torus^{I'}$.
If, for convenience (and without loss of generality),
we number the coordinates so that $I'=\{1,\ldots,k\},$ then
$\Gamma$ has a basis consisting of the rows of the block matrix.
\begin{equation}
\label{eq-splitdirac}
\left(
\begin{array}
[c]{llll}
\Pi_{\Gamma,M} & 0 & 1 & -\beta^t\\
\beta & 1 & 0 & 0
\end{array}
\right).
\end{equation}
where the $1$'s represent identity matrices.

We may see from this description of $\Gamma$ that the crossed product
algebra $\cala_{\Gamma,M}$ is generated by $n$ unitary elements which
satisfy the commutation relations associated with the Poisson
structure
\begin{equation}
\label{eq-poisson}
\left(
\begin{array}
[c]{ll}
\Pi_{\Gamma,M} & -\beta^t\\
\beta & 0
\end{array}
\right).
\end{equation}

The Dirac structure corresponding to this Poisson structure has a
basis given by the rows of the block matrix
\begin{equation}
\label{eq-poissondirac}
\left(
\begin{array}
[c]{llll}
\Pi_{\Gamma,M} & -\beta^t & I & 0\\
\beta & 0 & 0 & I
\end{array}
\right),
\end{equation}
which is obtained from \eqref{eq-splitdirac} by
interchange of the second and fourth columns, i.e. by the action
of the $O(n,n|\integers)$ element $\sigma_I$.

In addition to reproving part of Lemma \ref{lem-comp}, we have thus
proven:

\begin{prop}
\label{prop-quant}
If $\Gamma$ is a constant Dirac structure on $\torus^n$, and if
$\reals^{I'}$ is a complement to $C(\Gamma)$, then the crossed product
algebra $\cala_{\Gamma,M}$ is isomorphic to the quantum torus algebra
$\cala_{\sigma_I(\Gamma)}$.
\end{prop}

\begin{rk}
{\em Although $\cala_{\Gamma,M}$ depends
  only on $\Gamma$ and $M$, the Poisson structure $\sigma_I(\Gamma)$
  also depends on the choice of a torus complementary to $M$.
 However, changing the complement
  changes the Poisson structure by one with integer coefficients,
  which does not affect the quantization.
}
\end{rk}

\begin{thm}
\label{thm-quantization}
 The Morita equivalence class of $A _{\Gamma,M}$\ is independent of
the choice of $M$.
\end{thm}
\pf
We give an outline of the proof, with analytic details to appear in
\cite{ta:lemma}.
It will be an instance of a quantized
version of the fact (see \cite{mu-re-wi:equivalence}) that Morita
equivalent groupoids have Morita equivalent $C^*$-algebras.

Let $M_i$ ($i=1,2$) be transversal tori, used in Definition
\ref{defi-algebra} to produce algebras $\cala_{\Gamma,M_i}$.
These algebras are quantized versions of the
groupoid algebras of the transformation groupoids
$G_i=M_i\times\integers^k\arrows M_i$, which are
the restrictions to $M_1$ and $M_2$ of the transformation groupoid
$G=\torus^n\times \reals^k\arrows \torus^n$ associated to the
characteristic foliation.  An equivalence of groupoids between the $G_i$
is
given by the ``bibundle''
$M_1\stackrel{\mu_2}{\leftarrow}Q\stackrel{\mu_2}{\rightarrow} M_2$,
where $Q$ is the set of morphisms in
$G$ with target in $M_1$ and source in $M_2$, the ``moment maps''
$\mu_i:Q\to M_i$ are given by the target and source maps of $G$, and
the (free) actions
of $G_1$ and $G_2$ on $Q$ are given by left and right
$G$-multiplication respectively.  It is clear that these satisfy the
equivalence condition that the orbits of each action are the fibres
of the moment map of the other.

If our groupoid algebras were not quantized, the corresponding bimodule
would simply be $C(Q)$, as in \cite{mu-re-wi:equivalence}.
  To adapt to the quantized
groupoid algebras, we must also quantize $C(Q)$.  We do this
by observing first that the $\mu_i$'s are
covering maps, so we can pull back the Poisson structure
$\Pi_{\Gamma,M_1}$ by $\mu_1$ to give a Poisson structure
$\Pi_Q$ on $Q$.  Since the Poisson structures on $M_1$ and
$M_2$ are obtained from a invariant transverse bivector to the
characteristic foliation, the source map $\mu_2$ is then anti-Poisson.

We now use $\Pi_Q$ to quantize the functions on $Q$, obtaining an
algebra $\cala_Q$ into which $\cala_{\Gamma,M_1}$ and the opposite
algebra to $\cala_{\Gamma,M_1}$ embed as mutual commutants, making
$\cala_Q$ into a Morita equivalence bimodule.  The Hilbert module
structure is a ``quantized'' version of the one used in
\cite{mu-re-wi:equivalence}.
\qed

Theorem \ref{thm-quantization} makes the
following definition a valid one.

\begin{defi}
If $\Gamma$ is a constant Dirac structure on a torus, the {\bf
  quantization} of $\Gamma$ is the Morita equivalence class of
  algebras containing $\cala_{\Gamma,M}$ for any subtorus $M$ which is
  a complete transversal to the characteristic foliation of $\Gamma$.
\end{defi}

\begin{ex}
{\em
If $\Gamma$ is already a Poisson structure, we must take $M$ to be a
point, and so we recover the usual quantization of Poisson
structures.   On the other hand, if $\Gamma=F\oplus F^\circ$, where $F$
is a (constant) foliation, then the transverse Poisson structure is
zero, and we obtain the usual foliation algebra associated to a
complete transversal.  When $n=2$, any odd Dirac structure is of this
form, and the quantized algebra is a ``rotation algebra.''  Then
$\sigma_{\{1\}}(\Gamma) $ is a Poisson
structure whose quantization is the same algebra.  This gives a
geometric analog of the equivalence between rotation algebras and
quantum 2-torus algebras.
}
\end{ex}

\begin{rk}
{\em
Although the translation group of $\torus^n$ preserves the constant
Dirac structure $\Gamma$, it does not act on the algebra
$\cala_{\Gamma,M}$.  In some sense which it would be interesting to
study, it acts ``up to Morita equivalence.''}
\end{rk}

\section{Invariance under $O(n,n\vert \zz)$ Action}
\label{sec-invariance}

We are now ready to prove our main theorem.  Given the
Morita equivalences established in Section \ref{sec-quantization},
the proof will be very short.

\begin{thm}
\label{thm-main}
If two Dirac structures on $\torus^n$ are in the same orbit of
the $O(n,n|\mathbb{Z})$\ action, then their quantizations are the same
Morita equivalence class.
\end{thm}

\pf
In view of Lemma \ref{lem-generators}, it suffices to show that the
quantization is unchanged under the action of $\rho(GL(n,\integers))$
and of $\sigma_{\{1\}}$.

Any $A\in GL(n,\integers)$ acts on $\torus ^n$, and the invariant
nature of our construction implies that the algebras
$\cala_{\Gamma,M}$ and $\cala_{\rho(A)(\Gamma),A(M)}$
are isomorphic.

For $\sigma_{\{1\}}$, our argument is based on the following diagram.

$$
\begin{diagram}
\node{\Gamma}\arrow{s,l,T}{\sigma _{I}}\node{ \sigma _{ \{ 1\}
}(\Gamma)} \arrow{w,t,T}{\sigma _{\{ 1\}
}}\arrow{sw,b,T,..,2}{\sigma _{I \cup \{ 1 \}}}
\\
\node{\sigma _{I}(\Gamma)}
\end{diagram}
$$

To prove that $\Gamma$ and $\sigma_{\{1\}}\Gamma$ have the same
quantization, we first apply Lemma \ref{lem-nullity} to conclude
that $N(\Gamma)$ and
$N(\sigma_{\{1\}}(\Gamma))$ differ by $\pm 1$, which allows us to assume
that
$N(\sigma_{\{1\}}(\Gamma))=N(\Gamma) + 1$.

We now apply Proposition
\ref{prop-quant} and represent the quantization of $\Gamma$ by
the quantum torus algebra $\cala_{\sigma_I(\Gamma)}$, where
$|I|=N(\Gamma)$ and $\sigma_I(\Gamma)$ is a Poisson
structure.  Now $\sigma _{I \cup \{ 1 \}}(\sigma_{\{1\}}(\Gamma)) =
  \sigma_I(\Gamma)$ is a Poisson structure, so by Lemma
  \ref{lem-nullity} and Corollary \ref{co-length}, $|I\cup
  \{1\}|=N(\sigma_{\{1\}}(\Gamma)), $ and hence the quantization of
$\sigma_{\{1\}}(\Gamma)$ is also represented by
$\cala_{\sigma_I(\Gamma)}$.
\qed

In view of the $O(n,n|\integers)$ equivariance of the identification
of Poisson structures with their graphs, we have:

\begin{co}
If $\Pi\in {\calt}_n$ and $g\in SO(n,n|\integers)$ are such
that $g\cdot\Pi$ is defined, then the algebras $\cala_{\Pi}$ and
$\cala_{g\cdot\Pi}$ are Morita equivalent.
\end{co}

This is exactly extension of Theorem \ref{thm-risc}
 conjectured in \cite{ri-sc:morita} and proved
in \cite{li:strong}.

\section{Final remarks}
\label{sec-final}

\begin{rk}
{\em
In \cite{sc:morita}, A. Schwarz defines a refined notion of Morita
equivalence---complete Morita equivalence, which requires that there be a
connection on the Hilbert bimodule, compatible with
constant curvature connections on the noncommutative tori of both
sides. He showed that any two noncommutative tori must sit in the
same orbit of $SO(n,n|\integers)$ if they are completely Morita
equivalent. The statement of our main theorem \ref{thm-main} is also true
for complete Morita equivalence. The key point is that the  bimodule
constructed in \cite{ta:lemma} 
can be identified with a Heisenberg module as defined by
Rieffel in \cite{ri:high-mod}.  The complete Morita equivalence then
follows from Rieffel's result 
on Heisenberg modules . The relation between our bimodule
and the bimodule constructed in \cite{ri-sc:morita} will be studied in
\cite{ta:lemma}.
}
\end{rk}

\begin{rk}
{\em
Our proof of Theorem \ref{thm-main} is not optimal.
We would like to find  a uniform construction which takes
an  $O(n,n|\integers)$ element and produces a Morita equivalence
bimodule.  The
possibility of such a construction is related to the structure  of the
``Picard groupoid'' in which the objects are pairs $(\Gamma,M)$ and the 
morphisms are isomorphism
classes of Morita equivalence bimodules between the corresponding
algebras.  
}
\end{rk}

\begin{rk}
{\em Our approach in this paper may be related to open string theory
  in the following way.  In string theory, a subtorus $M$ may be
  considered as a ``$k$-brane.''  Strings ending on this brane are
  geodesics perpendicular to $M$ with respect to some background
  metric.  Therefore they lie in a foliation transverse to $M$.  If we
  are also given a ``$B$-field'' on $M$, represented by a constant
  bivector, then there is a unique Dirac structure $\Gamma$ on the
  torus for which the characteristic foliation $C(\Gamma)$ is the 
  given one and the transverse  Poisson structure $\Pi_{\Gamma,M}$ is
  the $B$-field.  The field
theory of of this open string is described (see, for example, 
\cite{ka:kronecker}) by the crossed product
algebra in Definition \ref{defi-algebra}.  

The action of $\sigma _{\{ i\}}$, can be seen as a t-duality in
the $i$-th direction. Lemma \ref{lem-generators} in this paper shows
that $O(n,n|\mathbb {Z})$ is the full group generated by
t-duality. And Theorem \ref{thm-main} shows that different open
string theories related by t-duality are Morita equivalent.}
\end{rk}

\begin{rk}
{\em There is another geometrical explanation of our Definition
\ref{defi-algebra}, which shows a relation to algebraic geometry and
homological mirror symmetry.

For a torus $\torus ^n$, we consider the moduli space of flat connections
on its trivial line bundle. Each such connection can be identified
with its holonomy,  a
homomorphism from the fundamental group of $\torus ^n$ to the
structure group $\torus ^1$.   It is easy to
see that the set of all flat connections is an $n$-torus,
but its typical tangent space is naturally identified with the
typical {\em co}tangent space of $\torus ^n$.  Therefore, we call this
new torus the 
{\bf dual torus of $\torus ^n$}, written as $\widehat {\torus ^n}$.

When a Dirac structure $\Gamma \in \cald _n$ has positive nullity, it
is not a Poisson structure.
However, when we fix a transversal $M$ to $C(\Gamma)$,
it is not difficult to see that $\Gamma$ defines a Poisson
structure on an n-dimensional subtorus $M\times M^0$ in $\torus
^n \times \widehat {\torus ^n}$. In particular, when $\Gamma$ is $\{0\}\times
\reals ^n$, it defines the zero Poisson structure on $\widehat {\torus
^n}$.

The group $O(n,n|\integers)$ acts on $\torus ^n \times \widehat {\torus
^n}$ naturally. The subgroup $\rho (GL(n, \integers))$ fixes the
product decomposition of $\torus ^n \times \widehat {\torus^n} $, while $\sigma
_{\{i\}}$ maps the $i$-th component torus to its dual.

The picture above can be seen as a real analog of the theory of
Mukai transformations \cite{mu:duality}
of abelian varieties.   We hope to elucidate
this connection and its
relation to Kontsevich's homological mirror symmetry
\cite{ko:homological}
in a future 
publication.
}
\end{rk}


\begin{thebibliography}{99}


\bibitem{bl-ge:quantization}
Block, J., and Getzler, E., Quantization of foliations, {\em
Proceedings of the XXth International Conference on
   Differential Geometric Methods in Theoretical Physics}, Vol. 1, 2
(New York, 1991), World Scientific, River Edge, NJ, 1992, 471--487.



\bibitem {C}
Connes, A., {\em Noncommutative Geometry}, Academic Press, 1994.

\bibitem {co-do-sc:nc-matrix}
Connes, A., Douglas, M.R., and Schwarz, A.,
Noncommutative Geometry and Matrix Theory: Compactification on
Tori, {\em J. High Energy Phys.}, 9802(1998)003


\bibitem{co:dirac}
Courant, T.J., Dirac manifolds, {\em Trans. A.M.S.} {\bf 319} (1990),
631--661.


\bibitem {el:k-theory}
Elliott, G. A., On the K-theory of the $C^{*}$ algebras generated by a
projective representation of a torsion-free
discrete
abelian group, {\em Operator Algebras and Group Representations}, 157-184,
Pitman, London, 1984.

\bibitem{ka:kronecker}
Kajiura, H., Kronecker foliation, $D1$-branes and Morita equivalence of
noncommutative two-tori, {\em J. High Energy Phys.} {\bf 8} (2002),
no.~50.

\bibitem{ko:homological}
Kontsevich, M.,  Homological algebra of mirror symmetry. 
{\em Proceedings of the International Congress of Mathematicians, 
Vol. 1, 2 (Z\"urich, 1994)} 120--139, Birkhäuser, Basel, 1995.


\bibitem{li:strong}
Li, H., Strong Morita equivalence of higher-dimensional noncommutative
tori, preprint math.OA/0303123.



\bibitem {li-sz:nc-du}
Lizzi, F., and Szabo, R., Noncommutative Geometry
and String Duality, {\em Corfu Summer Institute on Elementary Particle
Physics (Kerkyra, 1998, PRHEP-corfu98/073, 17 pp. (electronic)}, J. High
Energy Phys. Conf. Proc., J. High Energy Phys., Trieste, 1999.



\bibitem{mu-re-wi:equivalence}
Muhly, P.S., Renault, J.N., and Williams, D.P., Equivalence and
isomorphism for groupoid $C^*$-algebras, {\em J. Operator Theory} {\bf
17} (1987), 3--22.

\bibitem{mu:duality}
Mukai, S., Duality between $D(X)$ and $D(\hat X)$ with 
its application to Picard sheaves, {\em Nagoya Math. J.} {\bf 81} (1981),
  153--175.

\bibitem{ri:morita1}
Rieffel, M.A., Morita equivalence for
$C^{*}$-algebras and $W^{*}$-algebras, {\em J. Pure Appl. Algebra} {\bf 5}
(1974), 51-96.

\bibitem {ri:irrational}
Rieffel, M.A., {$C^{*}-$algebras associated with
irrational rotations}, {\em Pacific. J. Math.} {\bf 93} (1981), 415--429.


\bibitem {ri:high-mod}
Rieffel, M.A., Projective modules over higher-dimensional non-commutative
noncommutative tori, {\em Canadian J.Math.}40(1988), 257-338.

\bibitem{ri:deformation}
Rieffel, M.A., Deformation quantization of Heisenberg manifolds, {\em
Commun. Math. Phys.} {\bf 122} (1989), 531-562.



\bibitem{ri-sc:morita}
Rieffel, M.A., and Schwarz, A.,
Morita equivalence of multidimensional noncommutative tori, {\em
Int. J. Math.} {\bf 10} (1999), 289--299.

\bibitem{sc:morita}
Schwarz, A., Morita equivalence and duality, {\em
Lett. Math. Phys.} {\bf 50} (1999), 309--321.


\bibitem {ta:lemma}
Tang, X., Quantization of pseudo-\'etale groupoids (in preparation).


\bibitem{we:rotation}
Weinstein, A., Symplectic groupoids, geometric quantization,
and irrational rotation algebras, {\em Symplectic geometry, groupoids,
and integrable systems, S\'{e}minaire sud-Rhodanien
de g\'{e}om\'{e}trie \`{a} Berkeley (1989)}, P. Dazord and A.
Weinstein, eds., Springer-MSRI Series (1991), 281--290.



\bibitem{xu:noncommutative}
Xu, P., Noncommutative Poisson algebras, {\em Amer. J. Math.} {\bf
116} (1994), 101-125.

\end{thebibliography}
\end{document}